\documentclass[12pt]{amsart}
\usepackage{amsmath,amsfonts,amssymb,amsthm}
\usepackage{graphics,color,epsfig}

\newcommand{\R}{\mathbb R}
\newcommand{\C}{\mathbb C}
\newcommand{\Z}{\mathbb Z}

\newcommand{\cA}{\mathcal A}

\newcommand{\cF}{\mathcal F}

\newcommand{\eps}{\varepsilon}
\newcommand{\vhi}{\varphi}
\renewcommand{\Re}{\textnormal{Re~}}
\renewcommand{\Im}{\textnormal{Im~}}

\newtheorem{theorem}{Theorem}[section]
\newtheorem{definition}[theorem]{Definition}

\newtheorem{lemma}[theorem]{Lemma}
\newtheorem{proposition}[theorem]{Proposition}

\begin{document}
\title{Unique Ergodicity of Translation Flows}
\author{Yitwah Cheung and Alex Eskin}

\address{San Francisco State University \\
San Francisco, California}
\email{cheung@math.sfsu.edu}

\address{University of Chicago \\
Chicago, Illinois}
\email{eskin@math.uchicago.edu}
\subjclass{32G15, 30F30, 30F60, 37A25} 
\keywords{Abelian differentials, Teichm\"uller geodesics, 
Boshernitzan's criterion, Masur's criterion, slow divergence}
\date{\today}
\begin{abstract}
  This preliminary report contains a sketch of the proof of 
  the following result: a slowly divergent Teichm\"uller 
  geodesic satisfying a certain logarithmic law is determined 
  by a uniquely ergodic measured foliation.  
\end{abstract}
\maketitle

\section{Introduction}
Let $(X,\omega)$ be a translation surface.  
This means $X$ is a closed Riemann surface and $\omega$ 
  a holomorphic $1$-form on $X$.  
The line element $|\omega|$ induces a flat metric on $X$ which 
  has cone-type singularites at the zeroes of $\omega$ where 
  the cone angle is a integral multiple of $2\pi$.  
A \emph{saddle connection} in $X$ is a geodesic segment with 
  respect to the flat metric that joins a pair of zeroes of 
  $\omega$ without passing through one in its interior.  
Our main result is a new criterion for the unique ergodicity 
  of the vertical foliation $\cF_v$, defined by $\Re\omega=0$.  

\textbf{Teichm\"uller geodesics.}
The complex structure of $X$ is uniquely determined by the atlas 
  $\{(U_\alpha,\vhi_\alpha)\}$ of natural parameters away from 
  the zeroes of $\omega$ specified by $d\vhi_\alpha=\omega$.  
The evolution of $X$ under the Teichm\"uller flow is the family 
  of Riemann surfaces $X_t$ obtained by post-composing the charts 
  with the $\R$-linear map $z\to e^{t/2}\Re{z}+ie^{-t/2}\Im{z}$.  
It defines a unit-speed geodesic with respect to the Teichm\"uller 
  metric on the moduli space of compact Riemann surfaces normalised 
  so that Teichm\"uller disks have constant curvature $-1$.  
The Teichm\"uller map $f_t:X\to X_t$ takes rectangles to rectangles 
  of the same area, stretching in the horizontal direction and 
  contracting in the vertical direction by a factor of $e^{t/2}$.  
By a \emph{rectangle} in $X$ we mean a holomorphic map a product 
  of intervals in $\C$ such that $\omega$ pulls back to $dz$.  
All rectangles are assumed to have horizontal and vertical edges.  

Let $\ell(X_t)$ denote the length of the shortest saddle connection.  
Let $d(t)=-2\log\ell(X_t)$.  In this note we give a sketch of the
  proof of the following result.  

\begin{theorem}\label{thm:main}
There is an $\eps>0$ such that if $d(t)<\eps\log t+C$ for some 
  $C>0$ and for all $t>0$, then $\cF_v$ is uniquely ergodic.  
\end{theorem}

The hypothesis of Theorem~\ref{thm:main} can be formulated in 
  terms of the flat metric on $X$ without appealing to the 
  forward evolution of the surface.  
Let $h(\gamma)$ and $v(\gamma)$ denote the horizontal and vertical 
  components of a saddle connection $\gamma$, which are defined by 
$$h(\alpha)=\left|\Re\int_\gamma\omega\right|\quad\text{and}\quad
  v(\alpha)=\left|\Im\int_\gamma\omega\right|.$$
It is not hard to show that the following statements are equivalent.  
\begin{enumerate}
  \item[(a)] There is a $C>0$ such that for all $t>0$, $d(t)<\eps\log t+C.$ 
  \item[(b)] There is a $c>0$ such that for all $t>0$, $\ell(X_t)>c/t^{\eps/2}.$ 
  \item[(c)] There are constants $c'>0$ and $h_0>0$ such that for all 
        saddle connections $\gamma$ satisfying $h(\gamma)<h_0$, 
  \begin{equation}\label{strong:diophantine}
        h(\gamma)>\frac{c'}{v(\gamma)(\log v(\gamma))^\eps}.
  \end{equation}
\end{enumerate}


For any $p>1/2$ there are translation surfaces with nonergodic $\cF_v$ 
  whose Teichm\"uller geodesic $X_t$ satisfies the sublinear slow rate 
  of divergence $d(t)\le Ct^p$.  See \cite{Ch2}.  
Our main result asserts a logarithmic slow rate of divergence is enough 
  to ensure unique ergodicity of $\cF_v$.  

To illustrate our techniques, we give a new proof of the following 
  well-known result, which is a special case of Masur's criterion 
  for unique ergodicity \cite{Ma} which asserts that $\cF_v$ is 
  uniquely ergodic as soon as $X_t$ has an accumulation point in 
  the moduli space of compact Riemann surfaces.  
\begin{theorem}\label{thm:bosh}
If $\ell(X_t)$ does not approach zero as $t\to\infty$ 
  then $\cF_v$ is uniquely ergodic.  
\end{theorem}

\textbf{Remark.}
Theorem~\ref{thm:bosh} can also be proven using Boshernitzan's 
  criterion for interval exchange transformations, \cite{Ve1}, 
  \cite{Ve2}.\footnote{It should be pointed out that Masur's 
  result is misquoted in \cite{Ve2}.}
On the other hand, if $(X,\omega)$ is the suspension of a minimal 
  interval exchange transformation satisfying Boshernitzan's 
  criterion one can show that $\ell(X_t)$ does not tend to zero.  

It is not hard to give nonergodic examples of $\cF_v$ where the 
  ``$\liminf$'' of $d(t)$ goes to infinity arbitrarily slowly.  
\begin{theorem}\label{thm:liminf}
Given any function $r(t)\to\infty$ as $t\to\infty$ there exists 
  a Teichm\"uller geodesic $X_t$ with nonergodic $\cF_v$ and 
  $\liminf d(t)/r(t)<1$.  
\end{theorem}

\emph{Outline.} 
After introducing some notation and terminology in \S\ref{S:Generic} 
  we prove Theorem~\ref{thm:bosh} in \S\ref{S:Networks}.  
Then we explain what modifications are necessary to obtain a proof 
  of Theorem~\ref{thm:main}.  
A proof of Theorem~\ref{thm:liminf} is included in \S\ref{S:Liminf} 
  and can be read independent of the other sections.  


\section{Generic Points}\label{S:Generic}
If $\nu$ is a (normalised) ergodic invariant measure transverse to 
  the vertical foliation $\cF_v$ then for any horizontal arc $I$ 
  there is a full $\nu$-measure set of points $x\in X$ satisfying 
  \begin{equation}\label{eq:erg:ave}
    \lim \frac{^\#I\cap L_x}{|L_x|} = \nu(I) 
                            \quad\text{as}\quad |L_x|\to\infty
  \end{equation}
  where $L_x$ represents a vertical segment having $x$ as an endpoint.  
Given $I$, the set $E(I)$ of points satisfying (\ref{eq:erg:ave}) for 
  \emph{some} ergodic invariant $\nu$ has full Lebesgue measure.  
We refer to the elements of $E(I)$ as \emph{generic points} and 
  the limit in (\ref{eq:erg:ave}) as the \emph{ergodic average} 
  determined by $x$.  

\emph{Convention.} When passing to a subsequence $t_n\to\infty$ 
  along the Teichm\"uller geodesic $X_t$ we shall suppress the 
  double subscript notation and write $X_n$ instead of $X_{t_n}$.  
Similarly, we write $f_n$ instead of $f_{t_n}$.  

\begin{lemma}\label{lem:rectangle}
Let $x,y\in E(I)$ and suppose there is a sequence $t_n\to\infty$ 
  such that for every $n$ the images of $x$ and $y$ under $f_n$ 
  lie in a rectangle $R_n\subset X_n$ and the sequence of heights 
  $h_n$ satisfy $\lim h_ne^{t_n/2}=\infty$.  
Then $x$ and $y$ determine the same ergodic averages.  
\end{lemma}
\begin{proof}
One can reduce to the case where $f_n(x)$ and $f_n(y)$ lie at the 
  corners of $R_n$.  Let $n_-$ (resp. $n_+$) be the number of times 
  the left (resp. right) edge of $f_n^{-1}R_n$ intersects $I$.  
Observe that $n_-$ and $n_+$ differ by at most one so that since 
  $h_ne^{t_n/2}\to\infty$, the ergodic averages for $x$ and $y$ 
  approach the same limit.  
\end{proof}

Ergodic averages taken as $T\to\infty$ are determined by fixed 
  fraction of the tail: for any given $\lambda\in(0,1)$ 
  $$\frac{1}{T}\int_0^T f\to c \quad\text{ implies }\quad 
    \frac{1}{(1-\lambda)T}\int_{\lambda T}^T f\to c.$$  
This elementary observation is the motivation behind the following.  
\begin{definition}
A point $x$ is $K$-\emph{visible} from a rectangle $R$ if 
  the vertical distance from $x$ to $R$ is at most $K$ 
  times the height of $R$.  
\end{definition}

We have the following generalisation of Lemma~\ref{lem:rectangle}.  
\begin{lemma}\label{lem:visible}
If $x,y\in E(I)$, $t_n\to\infty$ and $K>0$ are such that for every $n$ 
  the images of $x$ and $y$ under $f_n$ are $K$-visible from some 
  rectangle whose height $h_n$ satisfies $h_ne^{t_n/2}\to\infty$, 
  then $x$ and $y$ determine the same ergodic averages.    
\end{lemma}

\section{Networks}\label{S:Networks}
To prove unique ergodicity we shall show that all ergodic averages 
  converge to the same limit.  
The ideas in this section were motivated by the proof of Theorem~1.1 
  in \cite{Ma}.  

\begin{definition}\label{def:reachable}
We say two points are $K$-\emph{reachable} from each other if 
  there is a rectangle $R$ from which both are $K$-visible.  
We also say two sets are $K$-\emph{reachable} from each other if 
  every point of one is $K$-reachable from every point of the other.  
\end{definition}

\begin{definition}\label{def:network}
Given a collection $\cA$ of subsets of $X$, we define an undirected 
  graph $\Gamma_K(\cA)$ whose vertex set is $\cA$ and whose edge 
  relation is given by $K$-reachability.  
A subset $Y\subset X$ is said to be $K$-\emph{fully covered} by $\cA$ 
  if every $y\in Y$ is $K$-reachable from some element of $\cA$.  
We say $\cA$ is a $K$-\emph{network} if $\Gamma_K(\cA)$ is connected 
  and $X$ is $K$-fully covered by $\cA$.  
\end{definition}

To prove Theorem~\ref{thm:bosh} we first need a proposition.  
\begin{proposition}\label{prop:network}
If $K>0$, $N>0$, $\delta>0$ and $t_n\to\infty$ are such that 
  for all $n$, there exists a $K$-network $\cA_n$ in $X_n$ 
  consisting of at most $N$ squares, each having measure at 
  least $\delta$, then $\cF_v$ is uniquely ergodic.  
\end{proposition}
\begin{proof}
Suppose $\cF_v$ is not uniquely ergodic.  Then we can find a 
  distinct pair of ergodic invariant measures $\nu_0$ and $\nu_1$ 
  and a horizontal arc $I$ such that $\nu_0(I)\neq\nu_1(I)$.  

We construct a finite set of generic points as follows.  
By allowing repetition, we may assume each $\cA_n$ contains exactly 
  $N$ squares, which shall be enumerated by $A(n,i), i=1,\dots,N$.  
Let $A_1\subset X$ be the set of points whose image under $f_n$ 
  belongs to $A(n,1)$ for infinitely many $n$.  
Note that $A_1$ has measure at least $\delta$ because it is a 
  descending intersection of sets of measure at least $\delta$.  
Hence, $A_1$ contains a generic point; call it $x_1$.  
By passing to a subsequence we can assume the image of $x_1$ 
  lies in $A(n,1)$ for all $n$.  
By a similar process we can find a generic point $x_2$ whose image 
  belongs to $A(n,2)$ for all $n$.  When passing to the subsequence, 
  the generic point $x_1$ retains the property that its image lies 
  in $A(n,1)$ for all $n$.  
Continuing in this manner, we obtain a finite set $F$ consisting of 
  $N$ generic points $x_i$ with the property that the image of $x_i$ 
  under $f_n$ belongs to $A(n,i)$ for all $n$ and $i$.  

Given a nonempty proper subset $F'\subset F$ we can always find 
  a pair of points $x\in F'$ and $y\in F\setminus F'$ such that 
  $f_n(x)$ and $f_n(y)$ are $K$-reachable from each other for 
  infinitely many $n$.  
This follows from the fact that $\Gamma_K(\cA_n)$ is connected.  
By Lemma~\ref{lem:visible}, the points $x$ and $y$ determine 
  the same ergodic averages for any horizontal arc $I$.  
Since $F$ is finite, the same holds for any pair of points in $F$.  

Now let $z_j$ be a generic point whose ergodic average is $\nu_j(I)$, 
  for $j=0,1$.  Since $X_n$ is $K$-fully covered by $\cA_n$, $z_j$ 
  will have the same ergodic average as some point in $F$, which 
  contradicts $\nu_0(I)\neq\nu_1(I)$.  
Therefore, $\cF_v$ must be uniquely ergodic.  
\end{proof}

\begin{proof}[Proof of Theorem~\ref{thm:bosh}]
Fix $\delta>0$ such that every saddle connection in $X_t$ has length 
  greater than $2\delta$ along some subsequence $t_n\to\infty$.  
Note that any immersed square of side at most $\delta$ is embedded 
  so that its Lebesgue measure is equal to its area.  
Consider the Delaunay triangulation of $X_n$.  (See \cite{MS}.)  
It has $2\nu$ triangles and $3\nu$ edges where $\nu$ depends only on 
  the genus and the number of singularities. 
Each triangle is contained in an immersed circumscribing disk of 
  radius at least $\delta$.  
It is not hard to see that every point in the disk is $1$-reachable 
  from the central square of side $\delta$ located at the circumcenter.  
Let $\cA_n=\{A_i\}_{i=1}^{2\nu}$ be the collection of central squares 
  and note that $X_n$ is $1$-fully covered by $\cA_n$.  

Let $\gamma_j$ be a Delaunay edge.  
It lies on the boundary of exactly two Delaunay triangles, each 
  having circumcenter on the perpendicular bisector of $\gamma_j$.  
Pick any square of side $\delta$ centered at some point on $\gamma_j$ 
  and divide it into $9$ nonoverlapping squares of side $\delta/3$.  
Two of the four corner squares are disjoint from $\gamma_j$.  
Let $A'_j$ (resp. $A''_j$) be the corner square that is disjoint from 
  $\gamma_j$ and lies on the same side of one of the circumcenters 
  (resp. the other circumcenter), possibly with $A''_j=A'_j$.  
One verifies easily that by adding these squares to the collection 
  $\cA_n$ we obtain a $1$-network of at most $8\nu$ squares, each 
  with area at least $\delta^2/9$.  
The theorem now follows from Proposition~\ref{prop:network}.  
\end{proof}

\section{Buffered Squares}
Our interest lies in the case where $\ell(X_t)\to0$ as $t\to\infty$.  
To prove of Theorem~\ref{thm:main} we shall need an analog of 
  Proposition~\ref{prop:network} that applies to a continuous family 
  of networks $\cA_t$ whose squares have dimensions going to zero.  
We also need to show that the slow rate of divergence gives us some 
  control on the rate at which the small squares approach zero.  

Recall that in the proof of Theorem~\ref{thm:bosh} it was essential 
  that the measure of the squares in the networks stay bounded away 
  from zero.  
This allowed us to find generic points that persist in the squares 
  of the networks along a subsequence $t_n\to\infty$.  
If the area of the squares tend to zero slowly enough, one can 
  still expect to find persistent generic points, with the help 
  of the following result from probability theory.  

\begin{lemma}\label{lem:PZ}
\textbf{(Paley-Zygmund \cite{PZ})}
If $A_n$ be a sequence of measureable subsets of a probability 
  space satisfying 
  \begin{enumerate}
    \item[(i)] $|A_n\cap A_m|\le K|A_n||A_m|$ for all $m>n$, and  
    \item[(ii)] $\sum |A_n| = \infty$ 
  \end{enumerate}
  then $|A_n \text{~i.o.}|\ge1/K$.  
\end{lemma}

\begin{definition}\label{def:buffer}
We say a rectangle is $\delta$-\emph{buffered} if it can be 
  extended in the vertical direction to a larger rectangle that 
  overlaps itself at most once and has area at least $\delta$.  
Here, we only require the product of the dimensions of the larger 
  rectangle be at least $\delta$, which obviously holds if its 
  Lebesgue measure exceeds $\delta$.  
\end{definition}

\begin{figure}
\begin{center}
\includegraphics{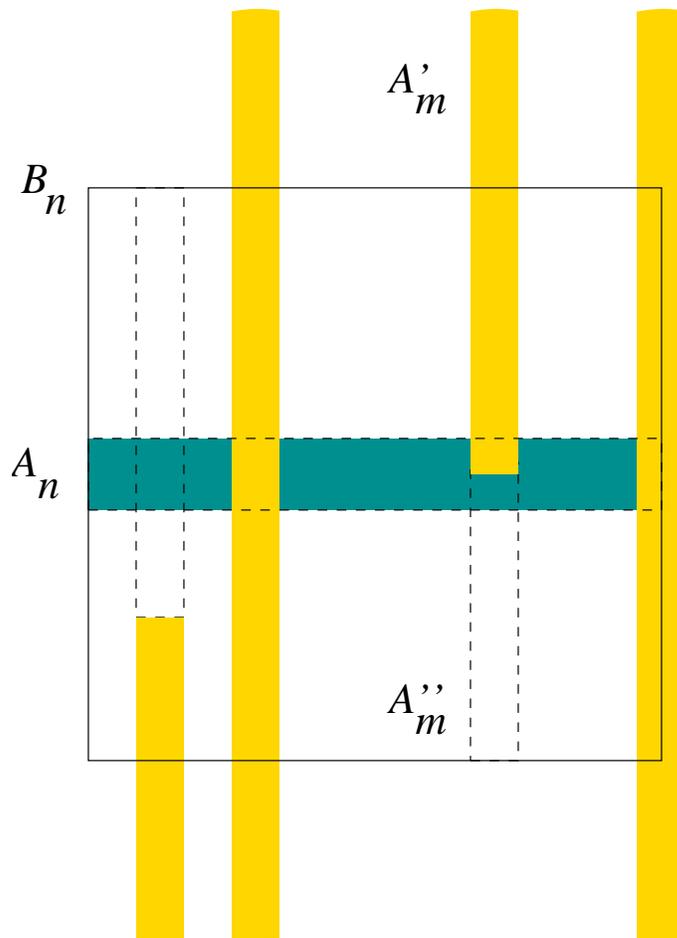}
\caption{For purposes of illustration, the rectangles are represented 
  by their images in $X_t$ where $t$ is the unique time when $B_n$ 
  maps to a square under the composition $f_t\circ f_n^{-1}$ of 
  Teichm\"uller maps.}\label{fig:buffer}
\end{center}
\end{figure}

\begin{proposition}\label{prop:spacing}
Suppose that for every $t>0$ we have a $\delta$-buffered square $A_t$ 
  embedded in $X_t$ with side $\sigma_t>c/t^{\eps/2}$, $0<\eps\le1$.  
Then there exists $t_n\to\infty$ such that conditions (i) and (ii) 
  of Lemma~\ref{lem:PZ} are satisfied by the sequence $f_n^{-1}A_n$ 
  of rectangles in $X$ for some $K$ depending only on $c$ and $\delta$.  
\end{proposition}
\begin{proof}
Let $(t_n)$ be any sequence satisfying the recurrence relation 
  $$t_{n+1}=t_n+\eps\log(t_{n+1}) \qquad t_0>1.$$  
Note that the function $y=y(x)=x-\eps\log x$ is increasing for $x>\eps$ 
  and has inverse $x=x(y)$ is increasing for $y>1$, from which it follows 
  that $(t_n)$ is increasing.  We have  $$\sigma_m e^{t_m-t_n} 
   >(ct_m^{-\eps/2})t_{n+1}^{\eps}\cdots t_m^{\eps/2}>ct_n^{\eps/2}.$$  
Let $B_n\supset A_n$ be the rectangle in $X_n$ that has the same 
  width as $A_n$ and area at least $\delta$.  
Let $\alpha(B_n)$ denote the product of the dimensions of $B_n$.  
Since $B_n$ overlaps itself at most once, $\alpha(B_n)\le2|B_n|\le2$.  
Therefore, the height of $B_n$ is $<2/\sigma_n<(2/c)t_n^{\eps/2}$, 
  which is less than $2/c^2$ times the height of the rectangle 
  $A'_m=f_n\circ f_m^{-1}A_m$, by the choice of $t_n$.  
Let $A''_m$ be the smallest rectangle containing $A'_m$ that has 
  horizontal edges disjoint from the interior of $B_n$.  
Its height is at most $1+4/c^2$ times that of $A_m$.  
For each component $I$ of $A_n\cap A''_m$ there is a corresponding 
  component $J$ of $B_n\cap A''_m$ (see Figure~\ref{fig:buffer}) 
  so that 
$$|A_n\cap A_m''|=\sum|I|\le\frac{\alpha(A_n)}{\alpha(B_n)}\sum|J|
  \le\delta^{-1}|A_n|\alpha(A''_m)<\frac{4+c^2}{\delta c^2}|A_n||A_m|$$ 
  giving (i).  
Choose $t_0$ large enough so that $t_{n+1}<2t_n$ for all $n$ and suppose 
  that for some $C>0$ and $n>1$ we have $t_n<Cn\log n\log\log n$.  Then 
\begin{align*}
  t_{n+1} &< t_n+\log t_{n+1} < t_n + \log t_n + \log 2 \\
     &< Cn\log n \log\log n + \log n + \log\log n + \log\log\log n + \log 2C \\
     &< Cn\log n \log\log n + \log n \log\log n \qquad\text{for $n\gg1$} \\
     &< C(n+1)\log(n+1)\log\log(n+1) 
\end{align*}
  so that $t_n\in O(n\log n\log\log n)$.  
Since $$\sum|A_n|=\sum\sigma_n^2>\sum\frac{c^2}{t_n^\eps}$$ 
  and $\eps\le1$, (ii) follows.  
\end{proof}

\begin{figure}
\begin{center}
\includegraphics{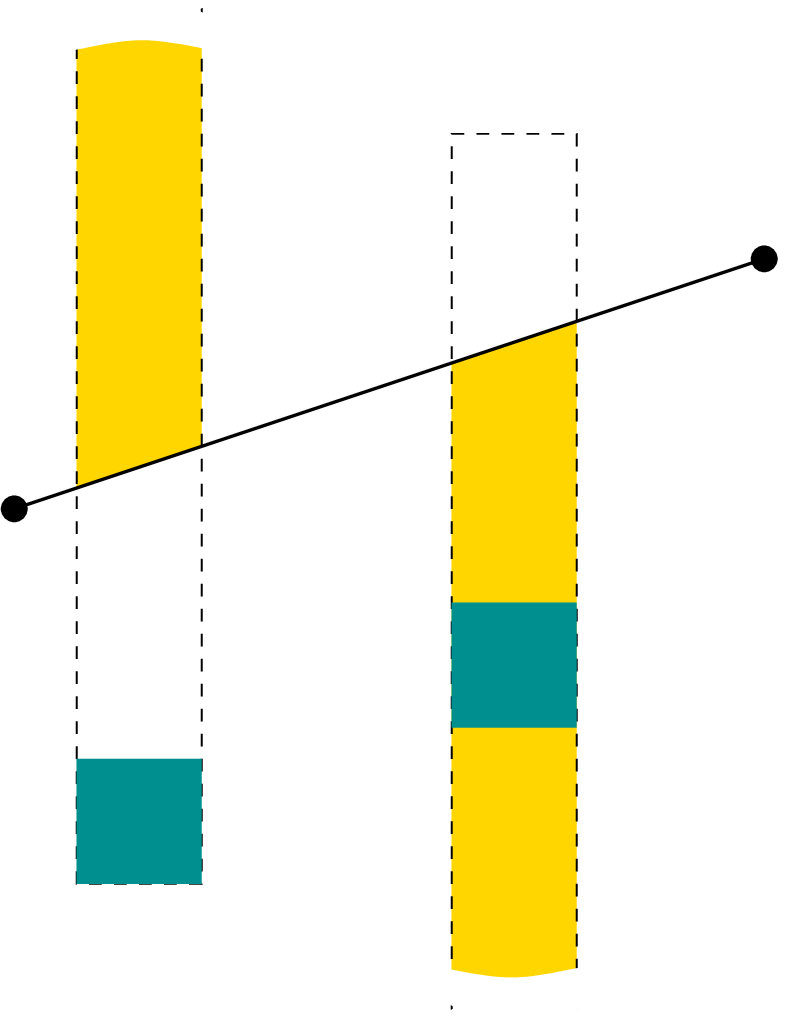}
\caption{Any rectangle contained the strip with largest area serves 
  as a buffer for any square of the same width contained in it.}
\label{fig:strip}
\end{center}
\end{figure}

There is a simple procedure for constructing buffered squares.  
Take any saddle connection $\gamma$ and decompose the surface 
  into a finite number of vertical strips, each of which is a 
  parallelogram with a zero on each of its vertical edges and 
  the remaining pair of edges contained in $\gamma$.  
The number of vertical strips is $\nu=2g-1+r$ where $g$ is the 
  genus of $X$ and $r$ the number of zeroes of $\omega$.  
Since one of the strips has area $\ge1/\nu$, any rectangle 
  containing the strip serves as a $1/\nu$-buffer for any 
  square of the same width contained in it.  
See Figure~\ref{fig:strip}.  

The condition (\ref{strong:diophantine}) prevents the slopes 
  of saddle connections from being too close to vertical.  
This allows for some control on the widths of vertical strips.  
\begin{proposition}\label{prop:strip}
There is a constant $c>0$ depending only on the constants of 
  the logarithmic rate such that for any $t>1$ and any saddle 
  connection $\gamma$ in $X_t$ of length at most $1$ there is 
  a vertical strip with a pair of edges contained in $\gamma$, 
  having area $\ge1/\nu$ and width $\ge c/t^{\nu\eps}$.  
\end{proposition}

By an argument similar to the proof of Theorem~\ref{thm:bosh} 
  one constructs a network of small buffered squares in $X_t$ 
  whose areas tend to zero no faster than $1/t^{2\nu\eps}$ 
  where $\eps$ is the coefficient of the logarithmic slow rate.  
Proposition~\ref{prop:strip} is used repeatedly on the edges 
  of the Delaunay triangulation to construct the individual 
  buffered squares of the network.  
(The main technical obstruction is that some of the edges of 
  Delaunay triangulation may be too long so that we cannot 
  apply Proposition~\ref{prop:strip}; however, long edges are 
  well understood--they have to cross long thin cylinders.)  
Invoking Proposition~\ref{prop:spacing} we pass to a subsequence 
  where the hypotheses of the lemma of Paley-Zygmund hold.  
An argument similar to the proof of Proposition~\ref{prop:network} 
  leads to the conclusion that $\cF_v$ is uniquely ergodic.  
Indeed, by the argument outlined above one can take $\eps=1/(2\nu)$ 
  in Theorem~\ref{thm:main}.\\

A complete proof of Theorem~\ref{thm:main} will appear elsewhere.

\section{Arbitrarily slow liminf divergence}\label{S:Liminf}
The examples with nonergodic vertical foliation are exhibited in 
  a well-known family of branched double covers of tori.  
For background regarding the construction, we refer the reader to 
  \cite{MT}.  

\begin{proof}[Proof of Theorem~\ref{thm:liminf}]
Let $X$ be the connected sum of $T^2=(\C/(\Z+i\Z),dz)$ with itself 
  along a horizontal arc $I$ of irrational length $\lambda\in(0,1)$.  
Let $W=\{\lambda+2m+2ni:m,n\in\Z,n>0\}$ and note that each element 
  in $W$ is a vector with irrational slope.  
The holonomy of a saddle connection $\gamma$ in $X$ joining the two 
  endpoints of the slit $I$ is of the form $\lambda+m+ni$ for some 
  $m,n\in\Z$ and since $\lambda$ is irrational, every such number is 
  the holonomy of some saddle connection in $X$ provided $n\neq0$.  
It is convenient to think of elements in $W$ as saddle connections.  
By concatenating $\gamma$ with its image under the hyperelliptic 
  involution, which interchanges the endpoints of $I$, we obtain 
  a simple closed curve $\Gamma$.  
The complement of $\Gamma$ has two connected components if and only 
  if $m$ and $n$ are both even integers; thus, every element of $W$ 
  determines a partition of $X$ into two sets of equal measure.  
For more details, see \cite{Ch1}.  

Fix a summable series of positive terms $\sum a_j<\infty$.  
Let $w_0\in W$.  Given $w_j\in W$ we can find an arbitrarily long 
  $v_j\in\Z+i\Z$ such that $|w_j\times v_j|<a_j$.  
(Here, $|z\times w|:=\Im(\Bar{z}w)$.)  
We set $w_{j+1}=w_j+2v_j$ and note that $w_{j+1}\in W$ and that 
  $|w_j\times w_{j+1}|<2a_j$.  
The symmetric difference between the partitions determined by $w_j$ 
  and $w_{j+1}$ is a union of parallelograms whose total area is 
  bounded above by the cross product.  
Applying Masur-Smillie's criterion for non-ergodicity \cite{MS}, 
  we see that the directions of $w_j$ converge to a nonergodic 
  direction $\theta$.  

Rotate so that $\theta$ is vertical and let $X_t$ be the evolution 
  of $X$ under the Teichm\"uller flow.  
Let $w_j^t$ be the image of $w_j$ in $X_t$.  
Let $t_j$ be the unique time $t=t_j$ when $w_j^t$ and $w_{j+1}^t$ 
  have the same length.  By choosing $v_j$ so that the vertical 
  (resp. horizontal) component of $w_{j+1}$ is at least twice 
  (resp. at most half) that of $w_j$ we can ensure that the angle 
  between $w_j^t$ and $w_{j+1}^t$ is at least $\cos^{-1}(4/5)$, 
  the angle between the vectors $1+2i$ and $2+i$.  
Hence, $\ell(X_t)^2\le|w_j^t||w_{j+1}^t|\le(5/3)|w_j\times w_{j+1}|$ 
  so that $d(t_j)\le-\log|w_j\times w_{j+1}|+C$.  
By choosing $v_j$ long enough we can make $d(t_j)/r(t_j)<1/2$, 
  which completes the proof.  
\end{proof}

\end{document}